\documentclass[12pt,openany]{article}
\usepackage{amsmath,amsthm,amsfonts,amssymb,amscd}
\usepackage{graphics}
\usepackage[utf8]{inputenc}

\renewcommand{\thefootnote}

\begin{document}

\centerline{\Large\bf On semiabelian categories in Functional Analysis and}

\medskip

\centerline{\Large\bf Topological Algebra}

\vglue .4in

\centerline{\bf Dinamérico P. Pombo Jr.}

\vglue .2in

\centerline{Instituto de Matemática}
\centerline{Universidade Federal do Rio de Janeiro}
\centerline{Caixa Postal 68530}
\centerline{21945-970 \,\, Rio de Janeiro, RJ\,\, Brasil}
\centerline{dpombo@terra.com.br}

\vglue .4in

\noindent{\bf Abstract.}\, In this work we discuss an elementary
self-contained presentation of the notion of a semiabelian category,
introduced by Ra$\overset{\vee}{\,\!\rm \i}\!$kov and Palamodov.
Fundamental examples of (non-abelian) semiabelian categories occuring in
Functional Analysis and Topological Algebra are treated in detail
throughout the paper.

\footnote{{\bf 2010 Mathematics Subject Classification}.\, 18-02, 18E05, 46H25, 22A05, 46A03, 46S10, 46A17.}
\footnote{{\bf Key words and phrases}.\, Preadditive categories, semiabelian categories, Functional Analysis, Topological Algebra.}

\vglue .3in

\noindent{\Large\bf \S1.\, Introduction}

\vglue .1in

Although abelian categories offer an adequate framework for dealing with
relevant aspects of Algebraic Geometry and Homological Algebra,
they do not contain many important categories which arise in
Functional Analysis and Topological Algebra, such as the categories
of Banach spaces, locally convex spaces, abelian Hausdorff topological
groups, topological modules over a topological ring and linearly topologized
modules over a topological ring. In [26], Ra$\overset{\vee}{\,\!\rm \i}\!$kov
introduced the so-called semiabelian categories (which were extensively
studied in [31]), which contain all abelian categories as well as certain
categories occuring in Functional Analysis and Topological Algebra
which are not abelian. Subsequently, Palamodov considered the notion
of a semiabelian category in a more general form (see [3,28,29,36]
for justifications of this assertion) and studied its application to certain
topics of the theory of locally convex spaces.

In this paper we discuss the concept of a semiabelian category in the sense
proposed by Palamodov. In order to make the exposition accessible to
non-specialists in Homological Algebra, we recall in the second section
a few basic notions and facts concerning general categories, the notion
of a category being a creation of Eilenberg and Mac Lane [11].
At this point we have avoided the customary use of the axiom of choice
in defining subobjects and quocients, which are important for our purposes,
following an approach due to Gabriel. The third section is devoted to basic
facts about preadditive categories which are needed for the comprehension
of the last section. Finally, in the last section, we use the background
presented in the preceding ones to introduce semiabelian categories.
After establishing some preparatory results, we concentrate our efforts
in proving a general criterion of bijectivity and an isomorphism theorem
in the context of semiabelian categories, in whose statements the concept
of a strict morphism plays a central role and which readily imply the known
isomorphism theorems for abelian categories. It should be mentioned that
various examples of relevant semiabelian categories which are not abelian,
encountered in Functional Analysis and Topological Algebra, are worked out
in detail throughout the text.

We would like to thank Yaroslav Kopylov for helpful comments which have
contributed to improve the presentation of the paper.

\vglue .3in

\noindent{\Large\bf \S2.\, Preliminaries on categories and examples}

\vglue .1in

For each \textit{category}\, $\mathcal C$\,[12,14], $\rm{Ob}(\mathcal C)$
will denote the class of \textit{objects\/} in $\mathcal C$ and,
for $A,B \in \rm{Ob}(\mathcal C)$, $\rm{Mor}_{\mathcal C}(A,B)$ will denote
the set of \textit{morphisms\/} from $A$ into $B$ in $\mathcal C$.
An element $u$ of $\rm{Mor}_{\mathcal C}(A,B)$ will be represented by
$u\colon A \to B$ or $A \overset{u}{\longrightarrow} B$.
For $A \in \rm{Ob}(\mathcal C)$, the identity morphism of $A$ will be
represented by $1_A$\,.

\medskip

Let us mention some examples of categories.

\bigskip

\noindent\textbf{Example 2.1.} The category Set whose objects are the sets
where, for $A,B \in \rm{Ob}(\rm Set)$, $\rm{Mor}_{\rm Set}(A,B)$ is the set
of mappings from $A$ into $B$.

\bigskip

\noindent\textbf{Example 2.2.} The category Top whose objects are the
topological spaces where, for $A,B \in \rm{Ob}(\rm Top)$,
$\rm{Mor}_{\rm Top}(A,B)$ is the set of continuous mappings from $A$ into $B$.

\bigskip

\noindent\textbf{Example 2.3.} The category Grp whose objects are the groups
where, for $A,B \in \rm{Ob}(\rm Grp)$, $\rm{Mor}_{\rm Grp}(A,B)$ is the group
of group homomorphisms from $A$ into $B$.

\bigskip

\noindent\textbf{Example 2.4.} For each ring $R$ with a non-zero identity
element, we can consider the category ${\rm Mod}_R$ whose objects are the
unitary left $R$-modules [5,20] where, for $A,B \in {\rm Ob}({\rm Mod}_R)$,
${\rm Mor}_{{\rm Mod}_R}(A,B)$ is the additive group of $R$-linear mappings
from $A$ into $B$. In the special case where $R$ is a field $K$
(resp. $R$ is the ring $\mathbb Z$ of integers), ${\rm Mod}_R$ is the
category of vector spaces over $K$ (resp. the category of abelian groups).

\bigskip

\noindent\textbf{Example 2.5.} For each topological ring $R$ with a non-zero
identity element, we can consider the category ${\rm Topm}_R$ whose objects
are the unitary left topological $R$-modules [34] where, for
$A,B \in {\rm Ob}({\rm Topm}_R)$, ${\rm Mor}_{{\rm Topm}_R}(A,B)$
is the additive group of continuous $R$-linear mappings from $A$ into $B$.
In the special case where $R$ is a topological field $K$ (resp. $R$ is the
ring $\mathbb Z$ of integers endowed with the discrete topology),
${\rm Topm}_R$ is the category of topological vector spaces over $K$\,[34]
(resp. the category of abelian topological groups [4,34]).

\bigskip

\noindent\textbf{Example 2.6.} For each topological ring $R$ with a non-zero
identity element, we can consider the category ${\rm Ltm}_R$ whose objects
are the unitary left linearly topologized $R$-modules [23,34] where, for
$A,B \in {\rm Ob}({\rm Ltm}_R)$, ${\rm Mor}_{{\rm Ltm}_R}(A,B)$ is the
additive group of continuous $R$-linear mappings from $A$ into $B$.
In the special case where $R$ is a discrete field $K$, ${\rm Ltm}_R$ is the
category of linearly topologized spaces over $\mathbb K$ (in the current
literature [17,21], linearly topologized spaces are also assumed to be
Hausdorff spaces).

\bigskip

\noindent\textbf{Example 2.7.} The category Ahtg whose objects are the
abelian Hausdorff topological groups where, for $A,B \in {\rm Ob}({\rm Ahtg})$,
${\rm Mor}_{\rm Ahtg}(A,B)$ is the abelian group of continuous group
homomorphisms from $A$ into $B$.

\bigskip

\noindent\textbf{Example 2.8.} The category Ban whose objects are the (real
or complex) Banach spaces [17,27] where, for $A,B \in {\rm Ob}({\rm Ban})$,
${\rm Mor}_{\rm Ban}(A,B)$ is the (real or complex) vector space of
continuous linear mappings from $A$ into $B$.

\bigskip

\noindent\textbf{Example 2.9.} The category Lcs whose objects are the (real
or complex) locally convex spaces [17,27] where, for
$A,B \in {\rm Ob}({\rm Lcs})$, ${\rm Mor}_{\rm Lcs}(A,B)$ is the (real or
complex) vector space of continuous linear mappings from $A$ into $B$.

\bigskip

\noindent\textbf{Example 2.10.} For each complete non-Archimedean
non-trivially valued field $K$, we can consider the category ${\rm Lcs}_K$
whose objects are the locally $K$-convex spaces [30,33] where, for
$A,B \in {\rm Ob}({\rm Lcs}_K)$, ${\rm Mor}_{{\rm Lcs}_K}(A,B)$ is the vector
space over $K$ of continuous $K$-linear mappings from $A$ into $B$.

\bigskip

\noindent\textbf{Example 2.11.} For each topological ring $R$ with a non-zero
identity element, we can consider the category ${\rm Borm}_R$ whose objects
are the bornological $R$-modules [24] where, for
$A,B \in {\rm Ob}({\rm Borm}_R)$, ${\rm Mor}_{{\rm Borm}_R}(A,B)$ in the
additive group of bounded $R$-linear mappings from $A$ into $B$.
In the special case where $R$ is a topological field $K$,
${\rm Mor}_{{\rm Borm}_R}(A,B)$ is the category of bornological vector spaces
over $K$ (the case where $K$ is a complete non-trivially valued field was
considered in [32]).

\bigskip

\noindent\textbf{Example 2.12.} The category Cbvs whose objects are the (real
or complex) convex bornological vector spaces [16] where, for
$A,B \in {\rm Ob}({\rm Cbvs})$, ${\rm Mor}_{\rm Cbvs}(A,B)$ is the (real or
complex) vector space of bounded linear mappings from $A$ into $B$.

\bigskip

\noindent\textbf{Example 2.13.} For each complete non-Archimedean non-trivially valued field $K$, we can consider the category ${\rm Cbvs}_K$ whose objects are the $K$-convex bornological vector spaces [1] where, for $A,B \in {\rm Ob}({\rm Cbvs}_K)$, ${\rm Mor}_{{\rm Cbvs}_K}(A,B)$ is the vector space over $K$ of bounded $K$-linear mappings from $A$ into $B$.

\bigskip

\noindent\textbf{Definition 2.14.} Let $\mathcal C$ be a category. The \textit{dual category\/} of $\mathcal C$, denoted by $\mathcal{C}^o$, is defined as follows:\, (a)\, ${\rm Ob}(\mathcal{C}^o) = {\rm Ob}(\mathcal C)$;\,\, (b)\, for $A,B \in {\rm Ob}(\mathcal{C}^o)$, ${\rm Mor}_{\mathcal{C}^o}(A,B) = {\rm Mor}_{\mathcal C}(B,A)$. The composition of morphisms
$$
{\rm Mor}_{\mathcal{C}^o}(A,B) \times {\rm Mor}_{\mathcal{C}^o}(B,C) \longrightarrow {\rm Mor}_{\mathcal{C}^o}(A,C)
$$
in $\mathcal{C}^o$ is defined by the composition of morphisms 
$$
{\rm Mor}_{\mathcal C}(C,B) \times {\rm Mor}_{\mathcal C}(B,A) \longrightarrow {\rm Mor}_{\mathcal C}(C,A)
$$
given in $\mathcal C$. Clearly, $(\mathcal{C}^o)^o = \mathcal C$.

\pagebreak

\noindent\textbf{Definition 2.15.} Let $\mathcal C$ be a category and $u \in {\rm Mor}_{\mathcal C}(A,B)$.\, $u$ is said to be \textit{injective\/} (resp. \textit{surjective}) if, for all $C \in {\rm Ob}(\mathcal C)$, the mapping $v \in {\rm Mor}_{\mathcal C}(C,A) \mapsto uv \in {\rm Mor}_{\mathcal C}(C,B)$ (resp. $w \in {\rm Mor}_{\mathcal C}(B,C) \mapsto wu \in {\rm Mor}_{\mathcal C}(A,C)$) is injective; $u$ is said to be \textit{bijective\/} if it is injective and surjective; $u$ is said to be an \textit{isomorphism\/} if there exists a $v \in {\rm Mor}_{\mathcal C}(B,A)$ such that $uv = 1_B$ and $vu = 1_A$ ($v$ is necessarily unique and denoted by $u^{-1})$.

Two objects $D$, $E$ in $\mathcal C$ are said to be \textit{isomorphic\/} if there is an isomorphism $w\colon D \to E$ in $\mathcal C$.

It is easily seen that the following assertions hold:

\smallskip

\noindent (a)\, In order that $u \in {\rm Mor}_{\mathcal C}(A,B)$ be injective (resp. $u \in {\rm Mor}_{\mathcal C}(A,B)$ be surjective), it is necessary and sufficient that $u \in {\rm Mor}_{\mathcal{C}^o}(B,A)$ be surjective (resp. $u \in {\rm Mor}_{\mathcal{C}^o}(B,A)$ be injective).

\smallskip

\noindent (b)\, If $u \in {\rm Mor}_{\mathcal C}(A,B)$, $v \in {\rm Mor}_{\mathcal C}(B,C)$ and $vu$ is injective (resp. $vu$ is surjective), then $u$ is injective (resp. $v$ is surjective).

\medskip

As a consequence of (b), every isomorphism is bijective. The following example furnishes a bijective morphism which is not an isomorphism.

\bigskip

\noindent\textbf{Example 2.16.} Let $A$ be the additive group of real numbers endowed with the discrete topology and $B$ the additive group of real numbers endowed with the usual topology, and let $u\colon A \to B$ be given by $u(x)=x$ for $x \in A$. Then $u \in {\rm Mor}_{\rm Ahtg}(A,B)$ is bijective, but $u$ is not an isomorphism.

\bigskip

\noindent\textbf{Definition 2.17} [13]. Let $\mathcal C$ be a category and $A \in {\rm Ob}(\mathcal C)$. A \textit{subobject\/} of $A$ is a pair $(A',i')$, where $A' \in {\rm Ob}(\mathcal C)$ and $i' \in {\rm Mor}_{\mathcal C}(A',A)$ is injective. Given two subobjects $(A',i')$, $(A'',i'')$ of $A$,\, $(A',i')$ is said to be bigger than or equal to $(A'',i'')$ (written $A' \ge A'')$ if there exists a morphism $u\colon A'' \to A'$ in $\mathcal C$ making the diagram
$$
\begin{matrix}
A''  &\overset{u}{\longrightarrow}& A'\\
i'' \searrow & & \swarrow i'\\
&A&
\end{matrix}  
$$
commutative. In this case, $u$ is injective and unique. For subobjects $(A',i')$, $(A'',i'')$ and $(A''',i''')$ of $A$, we have that $A' \ge A'$ and $A' \ge A'''$ if $A' \ge A''$ and $A'' \ge A'''$. Moreover, if $A' \ge A''$ and $A'' \ge A'$, with $i'' = i'u$ and $i' = i''v$ as above, then $u$ and $v$ are isomorphisms and $u^{-1} = v$.

\medskip

A \textit{quotient\/} of $A$ is a subobject of $A$,\, $A$ regarded as an object in the dual category ${\mathcal C}^o$. Therefore a quotient of $A$ is a pair $(Q',p')$, where $Q' \in {\rm Ob}(\mathcal C)$ and $p' \in {\rm Mor}(A,Q')$ is surjective. If $(Q',p')$, $(Q'',p'')$ are quotients of $A$,\, $(Q',p')$ is bigger than or equal to $(Q'',p'')$ (written $Q' \ge Q''$) if there exists a morphism $w\colon Q' \to Q''$ in $\mathcal C$ making the diagram
$$
\begin{matrix}
Q'  &\overset{w}{\longrightarrow}& Q''\\
p' \nwarrow & & \nearrow p''\\
&A&
\end{matrix}  
$$
commutative, $w$ being surjective and unique.

\bigskip

\noindent\textbf{Definition 2.18.} Let $\mathcal C$ be a category and $A,B \in {\rm Ob}(\mathcal C)$. Let $D \in {\rm Ob}(\mathcal C)$,\, $p_A \in {\rm Mor}_{\mathcal C}(D,A)$ and $p_B \in {\rm Mor}_{\mathcal C}(D,B)$. If, for each $C \in {\rm Ob}(\mathcal C)$, the mapping
$$
u \in {\rm Mor}_{\mathcal C}(C,D) \longmapsto (p_Au, p_Bu) \in {\rm Mor}_{\mathcal C}(C,A) \times {\rm Mor}_{\mathcal C}(C,B)
$$
is bijective, $D$ is said to be a \textit{product\/} of $A$, $B$ through $p_A$, $p_B$\,.

\smallskip

Two products of $A,B$, if they do exist, are isomorphic. For this reason, we shall refer to the product of $A, B$, denoted by $A \times B$. We shall say that \textit{finite products exist\/} in $\mathcal C$ if, for all $A,B \in {\rm Ob}(\mathcal C)$,\, $A \times B$ exists.

\medskip

Finite products exist in all the categories we have mentioned, as we shall
now see (more generally, it can be proved [25] that products exist in all
the categories considered in Examples 2.1 - 2.7 and in Examples 2.9 - 2.13).

\bigskip

\noindent\textbf{Example 2.19.} Finite products exist in Set.

\smallskip

In fact, let $A,B \in {\rm Ob}({\rm Set})$ be arbitrary. If $A \times B$ is the product set and $p_A\colon A \times B \to A$ and $p_B\colon A \times B \to B$ are the projections, it is obvious that, for each $C \in {\rm Ob}({\rm Set})$, the mapping
$$
u \in {\rm Mor}_{\rm Set}(C,A\times B) \longmapsto \big(p_A\circ u, p_B\circ u\big) \in {\rm Mor}_{\rm Set}(C,A) \times {\rm Mor}_{\rm Set}(C,B)
$$
is bijective.

\bigskip

\noindent\textbf{Example 2.20.} Finite products exist in Top.

\smallskip

In fact, let $A,B \in {\rm Ob}({\rm Top})$ be arbitrary and let $A \times B$,
$p_A$, $p_B$ be as in Example 2.19. If we endow $A \times B$ with the product
topology, it follows from Example 2.19 and Proposition 4, p.28 of [6] that
$A \times B$ is the product of $A,B$ through the continuous mappings
$p_A$, $p_B$\,.

\bigskip

\noindent\textbf{Example 2.21.} Finite products exist in Grp.

\smallskip

In fact, let $A,B \in {\rm Ob}({\rm Grp})$ be arbitrary and let $A \times B$, $p_A$, $p_B$ be as in Example 2.19. If we consider $A \times B$ endowed with its group structure, it follows from Example 2.19 that $A \times B$ is the product of $A,B$ through the group homomorphisms $p_A$, $p_B$\,.

\smallskip

By arguing as in Example 2.21, we have:

\bigskip

\noindent\textbf{Example 2.22.} Finite products exist in ${\rm Mod}_R$\,.

\bigskip

\noindent\textbf{Example 2.23.} Finite products exist in ${\rm Topm}_R$\,.

\smallskip

In fact, let $A,B \in {\rm Ob}({\rm Topm}_R)$ be arbitrary and let
$A \times B$, $p_A$, $p_B$ be as in Example 2.19. If we consider $A \times B$
endowed with its $R$-module structure and with the product topology, then
$A \times B \in {\rm Ob}({\rm Topm}_R)$ by Corollary 12.6 of [34].
Therefore, in view of Examples 2.20 and 2.22, \, $A \times B$ is the product
of $A,B$ through the continuous $R$-linear mappings $p_A$, $p_B$\,.

\bigskip

\noindent\textbf{Example 2.24.} Finite products exist in ${\rm Ltm}_R$\,.

\smallskip

It suffices to argue as in Example 2.23, Corollary 12.6 of [34] being
replaced by a remark on p.323 of [34].

\smallskip

By arguing as in Example 2.23, we get:

\bigskip

\noindent\textbf{Example 2.25.} Finite products exist in Ahtg.

\bigskip

\noindent\textbf{Example 2.26.} Finite products exist in Ban.

\smallskip

In fact, let $A,B \in {\rm Ob}({\rm Ban})$ be arbitrary and let $A \times B$, $p_A$, $p_B$ be as in Example 2.19. If we consider $A \times B$ endowed with its vector space structure and with the product norm, then $A \times B \in {\rm Ob}({\rm Ban})$ and it is easily verified that $A \times B$ is the product of $A,B$ through the continuous linear mappings $p_A$, $p_B$\,.

\smallskip

By arguing as in Example 2.23, with the pertinent modifications, we get:

\bigskip

\noindent\textbf{Example 2.27.} Finite products exist in Lcs.

\bigskip

\noindent\textbf{Example 2.28.} Finite products exist in ${\rm Lcs}_K$\,.

\bigskip

\noindent\textbf{Example 2.29.} Finite products exist in ${\rm Borm}_R$\,.

\smallskip

In fact, let $A,B \in {\rm Ob}({\rm Borm}_R)$ be arbitrary and let
$A \times B$, $p_A$, $p_B$ be as in Example 2.19.
If we consider $A \times B$ endowed with its $R$-module structure
and with the product bornology, it follows from Corollary 3 of [24] that
$A \times B \in {\rm Ob}({\rm Borm}_R)$ and is the product of $A,B$
through the bounded $R$-linear mappings $p_A$, $p_B$\,.

\bigskip

\noindent\textbf{Example 2.30.} Finite products exist in Cbvs.

\smallskip

In fact, let $A,B \in {\rm Ob}({\rm Cbvs})$ and let $A \times B$, $p_A$, $p_B$
be as in Example 2.19. If we consider $A \times B$ endowed with its vector
space structure and with the product bornology, it follows from Remark (1),
p.31 of [16] that $A \times B \in {\rm Ob}({\rm Cbvs})$ and is the product of
$A,B$ through the bounded linear mappings $p_A$, $p_B$\,.

\smallskip

By arguing as in Example 2.30, we get:

\bigskip

\noindent\textbf{Example 2.31.} Finite products exist in ${\rm Cbvs}_K$.

\smallskip

The dual notion to that considered in Definition 2.18 reads:

\bigskip

\noindent\textbf{Definition 2.32.} Let $\mathcal C$ be a category and $A,B \in {\rm Ob}(\mathcal C)$. Let $E \in {\rm Ob}(\mathcal C)$, $q_A \in {\rm Mor}_{\mathcal C}(A,E)$ and $q_B \in {\rm Mor}_{\mathcal C}(B,E)$. If, for each $C \in {\rm Ob}(\mathcal C)$, the mapping 
$$
u \in {\rm Mor}_{\mathcal C}(E,C) \longmapsto (uq_A, uq_B) \in {\rm Mor}_{\mathcal C}(A,C) \times {\rm Mor}_{\mathcal C}(B,C)
$$
is bijective, $E$ is said to be a \textit{coproduct\/} of $A,B$ through $q_A$, $q_B$\,.

\smallskip

Two coproducts of $A,B$, if they do exist, are isomorphic. For this reason, we shall refer to the coproduct of $A,B$, denoted by $A \sqcup B$. We shall say that \textit{finite coproducts exist\/} in $\mathcal C$ if, for all $A,B \in {\rm Ob}(\mathcal C)$,\, $A \sqcup B$ exists.

\bigskip

\noindent\textbf{Example 2.33.} Finite coproducts exist in Set.

\smallskip

In fact, let $A,B \in {\rm Ob}({\rm Set})$ be arbitrary.
Let $A \sqcup B$ be the sum of the sets $A,B$ and $q_A\colon A \to A \sqcup B$,
$q_B\colon B \to A \sqcup B$ the canonical mappings ([10], p.12).
Then it is easily seen that, for each $C \in {\rm Ob}({\rm Set})$,
the mapping
$$
u \in {\rm Mor}_{\rm Set}(A \sqcup B,C) \longmapsto (u \circ q_A, u \circ q_B) \in {\rm Mor}_{\rm Set}(A,C) \times {\rm Mor}_{\rm Set}(B,C)
$$
is bijective.

\bigskip

\noindent\textbf{Example 2.34.} Finite coproducts exist in Top.

\smallskip

In fact, let $A,B \in {\rm Ob}({\rm Top})$ be arbitrary and let $A \sqcup B$,
$q_A$, $q_B$ be as in Example 2.33. If we consider $A \sqcup B$ endowed with
the final topology for the mappings $q_A$, $q_B$\,, it follows from
Example 2.33 and Proposition 6, p.31 of [6] that $A \sqcup B$ is the
coproduct of $A,B$ through the continuous mappings $q_A$, $q_B$\,.

\bigskip

\noindent\textbf{Example 2.35.} Finite coproducts exist in Grp.

\smallskip

In fact, let $A,B \in {\rm Ob}({\rm Grp})$ be arbitrary. Let $A \sqcup B$ be
the free product of $A$, $B$ and $q_A\colon A \to A \sqcup B$,
$q_B\colon B \to A \sqcup B$ the canonical group homomorphisms.
Then, by Proposition 12.3 of [20], $A \sqcup B$ is the coproduct of $A,B$
through $q_A$, $q_B$\,.

\smallskip

In the sequel we shall see that finite coproducts exist in the categories
${\rm Mod}_R$\,, ${\rm Topm}_R$\,, ${\rm Ltm}_R$\,, Ahtg, Ban, Lcs,
${\rm Lcs}_K$\,, ${\rm Borm}_R$\,, Cbvs, ${\rm Cbvs}_K$\,.
More generally, it can be proved [25] that coproducts exist in all the
categories considered in Examples 2.1 - 2.6 and in Examples 2.9 - 2.13.

\vglue .3in

\noindent{\Large\bf \S3.\, Preliminaries on preadditive categories and}

\medskip

\noindent\qquad\quad\,{\Large\bf examples}

\vglue .1in

This section, strongly based on [2] and [14], is devoted to known facts about preadditive categories which will be needed in the sequel.

\bigskip

\noindent\textbf{Definition 3.1.} A category $\mathcal C$ is said to be \textit{preadditive\/} if the following axioms are satisfied:

\smallskip

\noindent (a)\, for all $A,B \in {\rm Ob}(\mathcal C)$, ${\rm Mor}_{\mathcal C}(A,B)$ is endowed with an abelian group structure  (the identity element of ${\rm Mor}_{\mathcal C}(A,B)$ will be denoted by $0_{AB}$);

\smallskip

\noindent (b)\, for all $A,B,C \in {\rm Ob}(\mathcal C)$, the composition of morphisms
$$
{\rm Mor}_{\mathcal C}(A,B) \times {\rm Mor}_{\mathcal C}(B,C) \longrightarrow {\rm Mor}_{\mathcal C}(A,C)
$$
is a $\mathbb{Z}$-bilinear mapping;

\smallskip

\noindent (c)\, there exists an $A \in {\rm Ob}(\mathcal C)$ such that $1_A = 0_{AA}$\,.

\medskip

Axiom (c) does not follow from the other axioms of a preadditive category: for instance, consider the category of non-trivial abelian groups.

It is obvious that the categories \, Set, Top and Grp are not preadditive, and that the categories ${\rm Mod}_R$\,, ${\rm Topm}_R$\,, ${\rm Ltm}_R$\,, Ahtg, Ban, Lcs, ${\rm Lcs}_K$\,, ${\rm Borm}_R$\,, Cbvs and ${\rm Cbvs}_K$ are preadditive.

\bigskip

\noindent\textbf{Remark 3.2.} Let $\mathcal C$ be a preadditive category and $A,B \in {\rm Ob}(\mathcal C)$. It can be shown that $A \times B$ exists if and only if $A \sqcup B$ exists; in this case, $A \times B$ and $A \sqcup B$ are isomorphic. 

\smallskip

Consequently, by Examples 22-31, finite coproducts exist in the categories ${\rm Mod}_R$\,, ${\rm Topm}_R$\,, ${\rm Ltm}_R$\,, Ahtg, Ban, Lcs, ${\rm Lcs}_K$\,, ${\rm Borm}_R$\,, Cbvs and ${\rm Cbvs}_K$\,.

\bigskip

\noindent\textbf{Remark 3.3.} Let $\mathcal C$ be a preadditive category and $u \in {\rm Mor}_{\mathcal C}(A,B)$. For each $C \in {\rm Ob}(\mathcal C)$ let $u^*$ (resp. $u_*$) be the group homomorphism from ${\rm Mor}_{\mathcal C}(C,A)$ into ${\rm Mor}_{\mathcal C}(C,B)$ (resp. from ${\rm Mor}_{\mathcal C}(B,C)$ into ${\rm Mor}_{\mathcal C}(A,C)$) given by $u^*(v) = uv$ for $v \in {\rm Mor}_{\mathcal C}(C,A)$ (resp. $u_*(w) = wu$ for $w \in {\rm Mor}_{\mathcal C}(B,C)$). It is clear that $u$ is injective (resp. surjective) if and only if the sequence
$$
0 \longrightarrow {\rm Mor}_{\mathcal C}(C,A) \overset{u^*}{\longrightarrow} {\rm Mor}_{\mathcal C}(C,B)\,\, ({\rm resp.}\, 0 \longrightarrow {\rm Mor}_{\mathcal C}(B,C) \overset{u_*}{\longrightarrow} {\rm Mor}_{\mathcal C}(A,C))
$$
of group homomorphisms is exact for each $C \in {\rm Ob}(\mathcal C)$ ([20], p.15).

\bigskip

\noindent\textbf{Definition 3.4.} Let $\mathcal C$ be a preadditive category. An object $A$ in $\mathcal C$ as in Definition 3.1(c) is said to be a \textit{zero\/} of $\mathcal C$.

\medskip

For $A \in {\rm Ob}(\mathcal C)$, the following conditions are equivalent:

\smallskip

\noindent (a)\, $A$ is a zero of $\mathcal C$;

\noindent (b)\, ${\rm Mor}_{\mathcal C}(A,A) = \{0_{AA}\}$;

\noindent (c)\, ${\rm Mor}_{\mathcal C}(A,C) = \{0_{AC}\}$ for all $C \in {\rm Ob}(\mathcal C)$;

\noindent (d)\, ${\rm Mor}_{\mathcal C}(C,A) = \{0_{CA}\}$ for all $C \in {\rm Ob}(\mathcal C)$.

\smallskip

Consequently, between two zeros of $\mathcal C$ there is only one isomorphism, namely, the unique morphism between them. For this reason, we shall refer to the zero of $\mathcal C$, denoted by $0$.

\bigskip

\noindent\textbf{Proposition 3.5.} \textit{Let $\mathcal C$ be a preadditive category and $u \in {\rm Mor}_{\mathcal C}(A,B)$. For $L \in {\rm Ob}(\mathcal C)$ and $i \in {\rm Mor}_{\mathcal C}(L,A)$, the following conditions are equivalent:}

\smallskip

\noindent (a)\quad $0 \longrightarrow {\rm Mor}_{\mathcal C}(C,L) \overset{i^*}{\longrightarrow} {\rm Mor}_{\mathcal C}(C,A) \overset{u^*}{\longrightarrow} {\rm Mor}_{\mathcal C}(C,B)$

\smallskip

\noindent\textit{is an exact sequence of group homomorphisms for each $C \in {\rm Ob}(\mathcal C)$.}

\smallskip

\noindent (b)\,\,\, (1)\, $i$ \textit{is injective};

\noindent\qquad (2)\, $ui = 0_{LB}$\,;

\noindent\qquad (3)\, \textit{for all $C \in {\rm Ob}(\mathcal C)$ and for all $v \in {\rm Mor}_{\mathcal C}(C,A)$ such that $uv = 0_{CB}$\,, there exists a $w \in {\rm Mor}_{\mathcal C}(C,L)$ such that $v = iw$}. 
 
\smallskip

Therefore, if $(L,i)$ and $(L',i')$ satisfy condition (a) of Proposition 3.5, then the subobjects $(L,i)$ and $(L',i')$ of $A$ are such that $L \ge L'$ and $L' \ge L$, and so are isomorphic.

\bigskip

\noindent\textbf{Definition 3.6.} Let $\mathcal C$ be a preadditive category and $u \in {\rm Mor}_{\mathcal C}(A,B)$. A pair $({\rm Ker}(u),i)$, where ${\rm Ker}(u) \in {\rm Ob}(\mathcal C)$ and $i \in {\rm Mor}_{\mathcal C}({\rm Ker}(u),A)$, is said to be the \textit{kernel\/} of $u$ if the sequence
$$
0 \longrightarrow {\rm Mor}_{\mathcal C}(C,{\rm Ker}(u)) \overset{i^*}{\longrightarrow} {\rm Mor}_{\mathcal C}(C,A) \overset{u^*}{\longrightarrow} {\rm Mor}_{\mathcal C}(C,B)
$$
of group homomorphisms is exact for each $C \in {\rm Ob}(\mathcal C)$.

The subobject $({\rm Ker}(u),i)$ of $A$, if it exists, is essentially unique, as we have just observed.

\bigskip

\noindent\textbf{Example 3.7.} Let $u \in {\rm Mor}_{{\rm Mod}_R}(A,B)$ and $C \in {\rm Ob}({\rm Mod}_R)$ be arbitrary. Put ${\rm Ker}(u) = \{x \in A; u(x)=0\}$ and let $i\colon {\rm Ker}(u) \to A$ be the inclusion mapping. Obviously, ${\rm Ker}(u) \in {\rm Ob}({\rm Mod}_R)$ and $i \in {\rm Mor}_{{\rm Mod}_R}({\rm Ker}(u),A)$. Since $i$ is an injective morphism in ${\rm Mod}_R$\,, Remark 3.3 guarantees the exactness of the sequence
$$
0 \longrightarrow {\rm Mor}_{{\rm Mod}_R}(C,{\rm Ker}(u)) \overset{i^*}{\longrightarrow} {\rm Mor}_{{\rm Mod}_R}(C,A).
$$

Now, let us consider the sequence
$$
{\rm Mor}_{{\rm Mod}_R}(C,{\rm Ker}(u)) \overset{i^*}{\longrightarrow} {\rm Mor}_{{\rm Mod}_R}(C,A) \overset{u^*}{\longrightarrow} {\rm Mor}_{{\rm Mod}_R}(C,B).
$$

If $w \in {\rm Mor}_{{\rm Mod}_R}(C,{\rm Ker}(u))$,
$$
u^*(i^*(w)) = u^*(i\circ w) = u \circ (i \circ w) = (u \circ i)\circ w = 0_{{\rm Ker}(u)B} \circ w = 0_{CB}\,;
$$
thus $i^*(w) \in {\rm Ker}(u^*)$ and ${\rm Im}(i^*) \subset {\rm Ker}(u^*)$. On the other hand, if $v \in {\rm Ker}(u^*)$, $\{v(x);\linebreak x \in C\} \subset {\rm Ker}(u)$, and we can view $v$ as an $R$-linear mapping $w$ from  $C$ into ${\rm Ker}(u)$. Moreover, $i^*(w) = i \circ w = v$, and hence ${\rm Ker}(u^*) \subset {\rm Im}(i^*)$. Therefore ${\rm Im}(i^*) = {\rm Ker}(u^*)$, and we have proved that $({\rm Ker}(u),i)$ is the kernel of $u$.

\bigskip

\noindent\textbf{Example 3.8.} Let $u \in {\rm Mor}_{{\rm Topm}_R}(A,B)$
and $C \in {\rm Ob}({\rm Topm}_R)$ be arbitrary. Let ${\rm Ker}(u)$ and $i$
be as in Example 3.7 and consider ${\rm Ker}(u)$ endowed with the topology
induced by that of $A$. Then ${\rm Ker}(u) \in {\rm Ob}({\rm Topm}_R)$
([34], p.87) and $i \in {\rm Mor}_{{\rm Topm}_R}({\rm Ker}(u),A)$.
We claim that $({\rm Ker}(u),i)$ is the kernel of $u$.
Indeed, since $i$ is an injective morphism in ${\rm Topm}_R$\,,
Remark 3.3 guarantees the exactness of the sequence
$$
0 \longrightarrow {\rm Mor}_{{\rm Topm}_R}(C, {\rm Ker}(u)) \overset{i^*}{\longrightarrow} {\rm Mor}_{{\rm Topm}_R}(C,A).
$$

Now, let us consider the sequence
$$
{\rm Mor}_{{\rm Topm}_R}(C,{\rm Ker}(u)) \overset{i^*}{\longrightarrow}{\rm Mor}_{{\rm Topm}_R}(C,A) \overset{u^*}{\longrightarrow} {\rm Mor}_{{\rm Topm}_R}(C,B).
$$
As in Example 3.7, ${\rm Im}(i^*) \subset {\rm Ker}(u^*)$. On the other hand, if $v \in {\rm Ker}(u^*)$, consider $w$ as in Example 3.7. Then it is clear that $w \in {\rm Mor}_{{\rm Topm}_R}(C, {\rm Ker}(u))$ and $i^*(w)=v$; thus ${\rm Ker}(u^*) \subset {\rm Im}(i^*)$. Therefore ${\rm Im}(i^*) = {\rm Ker}(u^*)$ and our claim is justified.

\medskip

By arguing exactly as in Example 3.8 (see a remark on p.323 of [34]), one concludes:

\bigskip

\noindent\textbf{Example 3.9.} Every morphism in\, ${\rm Ltm}_R$\, has a kernel.

\bigskip

\noindent\textbf{Example 3.10.} Let $u \in {\rm Mor}_{{\rm Ahtg}}(A,B)$ be arbitrary and let ${\rm Ker}(u)$ and $i$ be as in Example 3.8. Then ${\rm Ker}(u) \in {\rm Ob}({\rm Ahtg})$, $i \in {\rm Mor}_{{\rm Ahtg}}({\rm Ker}(u),A)$ and, by arguing exactly as in Example 3.8, one shows that $({\rm Ker}(u),i)$ is the kernel of $u$.

\bigskip

\noindent\textbf{Example 3.11.} Let $u \in {\rm Mor}_{{\rm Ban}}(A,B)$ be arbitrary and let ${\rm Ker}(u)$ and $i$ be as in Example 3.8. Since ${\rm Ker}(u)$ is a closed subspace of $A$, ${\rm Ker}(u) \in {\rm Ob}({\rm Ban})$ under the norm induced by that of $A$. Finally, $i \in {\rm Mor}_{\rm Ban}({\rm Ker}(u),A)$ and, by arguing as in Example 3.8, one shows that $({\rm Ker}(u),i)$ is the kernel of $u$.

\medskip

By arguing exactly as in Example 3.8, one justifies the validity of the two assertions below:

\bigskip

\noindent\textbf{Example 3.12.} Every morphism in Lcs has a kernel.

\bigskip

\noindent\textbf{Example 3.13.} Every morphism in ${\rm Lcs}_K$ has a kernel.

\bigskip

\noindent\textbf{Example 3.14.} Let $u \in {\rm Mor}_{{\rm Borm}_R}(A,B)$ be
arbitrary and let ${\rm Ker}(u)$ and $i$ be as in Example 3.7. Consider
${\rm Ker}(u)$ endowed with the bornology induced by that of $A$;
then ${\rm Ker}(u) \in {\rm Ob}({\rm Borm}_R)$ ([24], Corollary 1.a)) and $i \in {\rm Mor}_{{\rm Borm}_R}({\rm Ker}(u),A)$. We claim that $({\rm Ker}(u),i)$ is the kernel of $u$. Indeed, let $C \in {\rm Ob}({\rm Borm}_R)$ be arbitrary. Since $i$ is an injective morphism in ${\rm Borm}_R$\,, Remark 3.3 furnishes the exactness of the sequence
$$
0 \longrightarrow {\rm Mor}_{{\rm Borm}_R}(C,{\rm Ker}(u)) \overset{i^*}{\longrightarrow} {\rm Mor}_{{\rm Borm}_R}(C,A).
$$

Now, let us consider the sequence
$$
{\rm Mor}_{{\rm Borm}_R}(C, {\rm Ker}(u)) \overset{i^*}{\longrightarrow} {\rm Mor}_{{\rm Borm}_R}(C,A) \overset{u^*}{\longrightarrow} {\rm Mor}_{{\rm Borm}_R}(C,B).
$$
As in Example 3.7, ${\rm Im}(i^*) \subset {\rm Ker}(u^*)$. On the other hand, if $v \in {\rm Ker}(u^*)$, take $w$ as in Example 3.7. Then it is clear that $w \in {\rm Mor}_{{\rm Borm}_R}(C,{\rm Ker}(u))$ and $i^*(w) = v$. Hence ${\rm Ker}(u^*) \subset {\rm Im}(i^*)$, and ${\rm Im}(i^*) = {\rm Ker}(u^*)$. Therefore $({\rm Ker}(u),i)$ is the kernel of $u$.

\medskip

By arguing exactly as in Example 3.14, one justifies the validity of the two assertions below:

\bigskip

\noindent\textbf{Example 3.15.} Every morphism in Cbvs has a kernel.

\bigskip

\noindent\textbf{Example 3.16.} Every morphism in ${\rm Cbvs}_K$ has a kernel.

\medskip

By duality (recall Proposition 3.5), the following result holds:

\bigskip

\noindent\textbf{Proposition 3.17.} \textit{Let $\mathcal C$ be a preadditive category and $u \in {\rm Mor}_{\mathcal C}(A,B)$. For $J \in {\rm Ob}(\mathcal C)$ and $j \in {\rm Mor}_{\mathcal C}(B,J)$, the following conditions are equivalent}:

\smallskip

\noindent (a)\quad $0 \longrightarrow {\rm Mor}_{\mathcal C}(J,C) \overset{j_*}{\longrightarrow} {\rm Mor}_{\mathcal C}(B,C) \overset{u_*}{\longrightarrow} {\rm Mor}_{\mathcal C}(A,C)$

\smallskip

\noindent\textit{is an exact sequence of group homomorphisms for each $C \in {\rm Ob}(\mathcal C)$.}

\smallskip

\noindent (b)\,\,\, (1)\, $j$ \textit{is surjective};

\noindent\qquad (2)\, $ju = 0_{AJ}$\,;

\noindent\qquad (3)\, \textit{for all $C \in {\rm Ob}(\mathcal C)$ and for all $w \in {\rm Mor}_{\mathcal C}(B,C)$ such that $wu = 0_{AC}$\,, there exists a $v \in {\rm Mor}_{\mathcal C}(J,C)$ such that $w = vj$}. 
 
\smallskip

Therefore, if $(J,j)$ and $(J',j')$ satisfy condition (a) of Proposition 3.17, then the quotients $(J,j)$ and $(J',j')$ of $B$ are such that $J \ge J'$ and $J' \ge J$, and so are isomorphic.

\bigskip

\noindent\textbf{Definition 3.18.} Let $\mathcal C$ be a preadditive category and $u \in {\rm Mor}_{\mathcal C}(A,B)$. A pair $({\rm Coker}(u),j)$, where ${\rm Coker}(u) \in {\rm Ob}(\mathcal C)$ and $j \in {\rm Mor}_{\mathcal C}(B, {\rm Coker}(u))$,  is said to be the \textit{cokernel\/} of $u$ if the sequence
$$
0 \longrightarrow {\rm Mor}_{\mathcal C}({\rm Coker}(u),C) \overset{j_*}{\longrightarrow} {\rm Mor}_{\mathcal C}(B,C) \overset{u_*}{\longrightarrow} {\rm Mor}_{\mathcal C}(A,C)
$$
of group homomorphisms is exact for each $C \in {\rm Ob}(\mathcal C)$.

\smallskip

The quotient $({\rm Coker}(u),j)$ of $B$, if it exists, is essentially unique, as we have just observed.

\bigskip

\noindent\textbf{Example 3.19.} Let $u \in {\rm Mor}_{{\rm Mod}_R}(A,B)$ be arbitrary and let $M$ be the submodule $\{u(x); x \in A\}$ of $B$. Put ${\rm Coker}(u) = B/M$ (regarded as an $R$-module) and let $\pi\colon B \to {\rm Coker}(u)$ be the canonical $R$-linear mapping. We claim that $({\rm Coker}(u),\pi)$ is the cokernel of $u$. Indeed, let $C \in {\rm Ob}({\rm Mod}_R)$ be arbitrary. Since $\pi$ is a surjective morphism in ${\rm Mod}_R$\,, the sequence
$$
0 \longrightarrow {\rm Mor}_{{\rm Mod}_R}({\rm Coker}(u),C) \overset{\pi_*}{\longrightarrow} {\rm Mor}_{{\rm Mod}_R}(B,C)
$$
is exact by Remark 3.3. Now, let us consider the sequence
$$
{\rm Mor}_{{\rm Mod}_R}({\rm Coker}(u),C) \overset{\pi_*}{\longrightarrow} {\rm Mor}_{{\rm Mod}_R}(B,C) \overset{u_*}{\longrightarrow} {\rm Mor}_{{\rm Mod}_R}(A,C).
$$

If $v \in {\rm Mor}_{\mathcal C}({\rm Coker}(u),C)$,
$$
u_*(\pi_*(v)) = u_*(v\circ\pi) = (v\circ\pi)\circ u = v \circ (\pi\circ u) = v \circ 0_{A {\rm Coker}(u)} = 0_{AC}\,;
$$
hence ${\rm Im}(\pi_*) \subset {\rm Ker}(u_*)$. On the other hand, let $w \in {\rm Ker}(u_*)$ and define $v\colon {\rm Coker}(u) \to C$ by $v(\pi(y)) = w(y)$ for $y \in B$ ($v$ is well defined since ${\rm Ker}(\pi) \subset {\rm Ker}(w))$. Then $v \in {\rm Mor}_{{\rm Mod}_R}({\rm Coker}(u),C)$ and $w = v \circ \pi = \pi_*(v)$; thus ${\rm Ker}(u_*) \subset {\rm Im}(\pi_*)$. Therefore  ${\rm Im}(\pi_*) = {\rm Ker}(u_*)$, and $({\rm Coker}(u),\pi)$ is the cokernel of $u$.

\bigskip

\noindent\textbf{Example 3.20.} Let $u \in {\rm Mor}_{{\rm Topm}_R}(A,B)$ be
arbitrary and let ${\rm Coker}(u)$ and $\pi$ be as in Example 3.19.
By considering ${\rm Coker}(u)$ endowed with the quotient topology,
it follows that ${\rm Coker}(u) \in {\rm Ob}({\rm Topm}_R)$
([34], Theorem 12.10) and $\pi \in {\rm Mor}_{{\rm Topm}_R}(B, {\rm Coker}(u))$. We claim that $({\rm Coker}(u),\pi)$ is the cokernel of $u$. Indeed, let $C \in {\rm Ob}({\rm Topm}_R)$ be arbitrary. Since $\pi$ is a surjective morphism in ${\rm Topm}_R$\,, the sequence
$$
0 \longrightarrow {\rm Mor}_{{\rm Topm}_R}({\rm Coker}(u),C) \overset{\pi_*}{\longrightarrow} {\rm Mor}_{{\rm Topm}_R}(B,C)
$$
is exact by Remark 3.3. Now, let us consider the sequence
$$
{\rm Mor}_{{\rm Topm}_R}({\rm Coker}(u),C) \overset{\pi_*}{\longrightarrow} {\rm Mor}_{{\rm Topm}_R}(B,C) \overset{u_*}{\longrightarrow} {\rm Mor}_{{\rm Topm}_R}(A,C).
$$
As in Example 3.19, ${\rm Im}(\pi_*) \subset {\rm Ker}(u_*)$. On the other hand, if $w \in {\rm Ker}(u_*)$, define $v$ as in Example 3.19. Then $v \in {\rm Mor}_{{\rm Topm}_R}({\rm Coker}(u),C)$. In fact, if $V$ is a neighborhood of $0$ in $C$ there is a neighborhood $U$ of $0$ in $B$ such that $w(U) \subset V$. Since $v(\pi(U)) = w(U)$, $\pi(U)$ being a neighborhood of $0$ in ${\rm Coker}(u)$ because $\pi$ is open, the continuity of $v$ is verified. Finally, since $w = \pi_*(v)$, the inclusion ${\rm Ker}(u_*) \subset {\rm Im}(\pi_*)$ is valid. Therefore ${\rm Im}(\pi_*) = {\rm Ker}(u_*)$, and our claim is proved.

\bigskip

\noindent\textbf{Example 3.21.} Let $u \in {\rm Mor}_{{\rm Ltm}_R}(A,B)$ be arbitrary and let ${\rm Coker}(u)$ and $\pi$ be as in Example 3.20. Then it is clear that ${\rm Coker}(u) \in {\rm Ob}({\rm Ltm}_R)$ and $\pi \in {\rm Mor}_{{\rm Ltm}_R}(B, {\rm Coker}(u))$. Moreover, by arguing exactly as in Example 3.20, one shows that $({\rm Coker}(u),\pi)$ is the cokernel of $u$.

\bigskip

\noindent\textbf{Example 3.22.} Let $u \in {\rm Mor}_{\rm Ahtg}(A,B)$ be arbitrary. Then the set $M = \{\overline{u(x); x \in A}\}$ is a closed subgroup of $B$, and hence the abelian group ${\rm Coker}(u) = B/M$, endowed with the quotient topology, is an object in Ahtg. If $\pi\colon B \to {\rm Coker}(u)$ denotes the canonical group homomorphism, which is a morphism in Ahtg, the argument used in Example 3.20 shows that $({\rm Coker}(u),\pi)$ is the cokernel of $u$.

\bigskip

\noindent\textbf{Example 3.23.} Let $u \in {\rm Mor}_{\rm Ban}(A,B)$ be arbitrary. Then the set $M = \{\overline{u(x); x \in A}\}$ is a closed subspace of $B$, and hence the vector space ${\rm Coker}(u) = B/M$, endowed with the quotient norm, is an object in Ban. If $\pi\colon B \to {\rm Coker}(u)$ denotes the canonical linear mapping, which is a morphism in Ban, the argument used in Example 3.20 shows that $({\rm Coker}(u),\pi)$ is the cokernel of $u$.

\medskip

By arguing as in Example 3.20, one shows the validity of the two assertions below:

\bigskip

\noindent\textbf{Example 3.24.} Every morphism in Lcs admits a cokernel.

\bigskip

\noindent\textbf{Example 3.25.} Every morphism in ${\rm Lcs}_K$ admits a cokernel.

\bigskip

\noindent\textbf{Example 3.26.} Let $u \in {\rm Mor}_{{\rm Borm}_R}(A,B)$ be arbitrary and let ${\rm Coker}(u)$ and $\pi$ be as in Example 3.19. By considering ${\rm Coker}(u)$ endowed with the quotient bornology, we have that ${\rm Coker}(u) \in {\rm Ob}({\rm Borm}_R)$ and $\pi \in {\rm Mor}_{{\rm Borm}_R}(B,{\rm Coker}(u))$. We claim that $({\rm Coker}(u),\pi)$ is the cokernel of $u$. Indeed, let $C \in {\rm Ob}({\rm Borm}_R)$ be arbitrary. Since $\pi$ is a surjective morphism in ${\rm Borm}_R$\,, Remark 3.3 implies the exactness of the sequence
$$
0 \longrightarrow {\rm Mor}_{{\rm Borm}_R}({\rm Coker}(u),C) \overset{\pi_*}{\longrightarrow} {\rm Mor}_{{\rm Borm}_R}(B,C).
$$

Now, let us show the exactness of the sequence
$$
{\rm Mor}_{{\rm Borm}_R}({\rm Coker}(u),C) \overset{\pi_*}{\longrightarrow} {\rm Mor}_{{\rm Borm}_R}(B,C) \overset{u_*}{\longrightarrow} {\rm Mor}_{{\rm Borm}R}(A,C).
$$
As in Example 3.19, ${\rm Im}(\pi_*) \subset {\rm Ker}(u_*)$. On the other hand, if $w \in {\rm Ker}(u_*)$, define $v$ as in Example 3.19. Then $v \in {\rm Mor}_{{\rm Borm}_R}({\rm Coker}(u),C)$ because $w$ is bounded and $v(\pi(L)) = w(L)$ for every bounded subset $L$ of $B$. Finally, since $w = \pi_*(v)$, the inclusion ${\rm Ker}(u_*) \subset {\rm Im}(\pi_*)$ is valid. Therefore ${\rm Im}(\pi_*) = {\rm Ker}(u_*)$, and our claim is proved.

\medskip

By arguing as in the preceding example, one justifies the validity of the two assertions below:

\bigskip

\noindent\textbf{Example 3.27.} Every morphism in Cbvs has a cokernel.

\bigskip

\noindent\textbf{Example 3.28.} Every morphism in ${\rm Cbvs}_K$ has a cokernel.

\bigskip

\noindent\textbf{Definition 3.29.} Let $\mathcal C$ be a preadditive category satisfying the following axiom:

\smallskip

(d)\, every morphism in $\mathcal C$ admits kernel and cokernel.

\smallskip

\noindent For $u \in {\rm Mor}_{\mathcal C}(A,B)$, we define the \textit{image\/} of $u$, denoted by ${\rm Im}(u)$, as ${\rm Im}(u) = {\rm Ker}({\rm Coker}(u))$, and the \textit{coimage\/} of $u$, denoted by ${\rm Coim}(u)$, as ${\rm Coim}(u) = {\rm Coker}({\rm Ker}(u))$
$$
\begin{matrix}
{\rm Ker}(u) \overset{i}{\longrightarrow} &A& \overset{j'}{\longrightarrow} {\rm Coim}(u)\\
&\,\,\downarrow u&\\
{\rm Im}(u) \underset{i'}{\longrightarrow} &B& \underset{j}{\longrightarrow} {\rm Coker}(u)
\end{matrix}
$$

As we have already seen, the preadditive categories ${\rm Mod}_R$\,, ${\rm Topm}_R$\,, ${\rm Ltm}_R$\,, Ahtg, Ban, Lcs, ${\rm Lcs}_K$\,, ${\rm Borm}_R$, Cbvs and ${\rm Cbvs}_K$ satisfy axiom (d).

\medskip

In the statement of the next result, the equality ${\rm Coim}(u)=A$ (resp. ${\rm Im}(u)=B$) will mean that the morphism $j'\colon A \to {\rm Coim}(u)$ (resp. $i'\colon {\rm Im}(u) \to B$) is an isomorphism.

\bigskip

\noindent\textbf{Proposition 3.30.} \textit{Let $\mathcal C$ be a category as in Definition $3.29$ and let $u \in {\rm Mor}_{\mathcal C}(A,B)$. Then the following assertions hold:}

\smallskip

\noindent (a)\, \textit{$u$ is injective if and only if \, ${\rm Ker}(u)=0$};

\smallskip

\noindent (b)\, \textit{$u$ is injective if and only if \, ${\rm Coim}(u)=A$};

\smallskip

\noindent (c)\, \textit{$u$ is surjective if and only if \, ${\rm Coker}(u)=0$};

\smallskip

\noindent (d)\, \textit{$u$ is surjective if and only if \, ${\rm Im}(u)=B$}.

\medskip

The following basic result, whose proof we shall include here, will be important for our purposes.

\noindent\textbf{Proposition 3.31.} \textit{Let $\mathcal C$ be a category as in Definition $3.29$ and let $u \in {\rm Mor}_{\mathcal C}(A,B)$. Then there exists a unique $\overline{u} \in {\rm Mor}_{\mathcal C}({\rm Coim}(u), {\rm Im}(u))$ such that $u = i'\,\overline{u}\,j'$, $i'$ and $j'$ being as in Definition $3.29$.}
$$
A \overset{j'}{\longrightarrow} {\rm Coim}(u) \overset{\overline{u}}{\longrightarrow} {\rm Im}(u) \overset{i'}{\longrightarrow} B.
$$

\medskip

\noindent\textbf{Proof.} In order to prove the uniqueness, assume that $u = i'\,\overline{u}\,j'$ and $u = i'\,\overline{u'}\,j'$, where $\overline{u},\overline{u'} \in {\rm Mor}_{\mathcal C}({\rm Coim}(u), {\rm Im}(u))$. Since $j'$ is surjective and $(i'\,\overline{u})j' = (i'\,\overline{u'})j'$, then $i'\,\overline{u} = i'\,\overline{u'}$. Thus $\overline{u} = \overline{u'}$ because $i'$ is injective.

\smallskip

Now, let us prove the existence. Indeed, let $({\rm Ker}(u),i)$ be the kernel of $u$ and $({\rm Coker}(u),j)$ the cokernel of $u$. Since $ui = 0_{{\rm Ker}(u)B}$ by Proposition 3.5, there is a morphism $v\colon {\rm Coim}(u) \to B$ in $\mathcal C$ such that $u = vj'$ (Proposition 3.17). On the other hand, since $ju = 0_{A\,{\rm Coker}(u)}$ by Proposition 3.17, then
$$
(jv)j' = j(vj') = ju = 0_{A\,{\rm Coker}(u)} = 0_{{\rm Coim}(u){\rm Coker}(u)} j'.
$$
Thus, by the surjectivity of $j'$, $jv = 0_{{\rm Coim}(u){\rm Coker}(u)}$. Therefore, by Proposition 3.5, there is a morphism $\overline{u}\colon {\rm Coim}(u) \to {\rm Im}(u)$ in $\mathcal C$ such that $v = i'\,\overline{u}$. Consequently, $u = vj' = i'\,\overline{u}\,j'$, which concludes the proof.

\medskip

The proofs of the next results may be found in [2].

\bigskip

\noindent\textbf{Proposition 3.32.} \textit{Let $\mathcal C$ be a category as in Definition $3.29$, $u \in {\rm Mor}_{\mathcal C}(A,B)$ and $v \in {\rm Mor}_{\mathcal C}(B,C)$. Then the following assertions hold:}

\smallskip

\noindent (a)\, \textit{${\rm Ker}(vu) \ge {\rm Ker}(u)$, and ${\rm Ker}(vu) \le {\rm Ker}(u)$ if $v$ is injective;}

\smallskip

\noindent (b)\, \textit{${\rm Coker}(vu) \ge {\rm Coker}(v)$, and ${\rm Coker}(vu) \le {\rm Coker}(v)$ if $u$ is surjective.}

\bigskip

\noindent\textbf{Corollary 3.33.} \textit{For $\mathcal C$, $u$ and $v$ as in Proposition $3.32$, the following assertions hold:}

\smallskip

\noindent (a)\, \textit{${\rm Coim}(vu) \le {\rm Coim}(u)$, and ${\rm Coim}(vu) \ge {\rm Coim}(u)$ if $v$ is injective;}

\smallskip

\noindent (b)\, \textit{${\rm Im}(vu) \le {\rm Im}(v)$, and ${\rm Im}(vu) \ge {\rm Im}(v)$ if $u$ is surjective.}

\bigskip

\noindent\textbf{Corollary 3.34.} \textit{Let $\mathcal C$ be a category as in Definition $3.29$ and $A \in {\rm Ob}(\mathcal C)$. Then the following assertions hold:}
\begin{itemize}
\item[{(a)}] \textit{if $(A',i')$ and $(A'',i'')$ are subobjects of $A$ such that $A' \le A''$, then ${\rm Coker}(i') \ge {\rm Coker}(i'')$;}
\item[{(b)}] \textit{if $(Q',p')$ and $(Q'',p'')$ are quotients of $A$ such that $Q' \le Q''$, then ${\rm Ker}(p') \ge {\rm Ker}(p'')$.}
\end{itemize}

\medskip

\noindent\textbf{Proof.} (a):\, Since $A' \le A''$, there is a morphism $u\colon A' \to A''$ in $\mathcal C$ such that $i' = i'' u$. Thus, by Proposition 3.32(b), ${\rm Coker}(i') \ge {\rm Coker}(i'')$.

\smallskip

\noindent (b):\, Follows from (a), by duality.

\vglue .3in

\noindent{\Large\bf \S4.\, Semiabelian categories}

\bigskip

\noindent\textbf{Definition 4.1}\, [22]. A category $\mathcal C$ is said to be \textit{semiabelian\/} if it is preadditive, satisfies axiom (d), as well as the following axioms:
\begin{itemize}
\item[{(e)}]\, finite products exist in $\mathcal C$;
\item[{(f)}]\, for all $A,B \in {\rm Ob}(\mathcal C)$ and for all $u \in {\rm Mor}_{\mathcal C}(A,B)$, the morphism $\overline{u}\colon {\rm Coim}(u) \to {\rm Im}(u)$ (Proposition 3.31) is bijective.
\end{itemize}

\bigskip

By Remark 3.2, finite coproducts exist in any semiabelian category.

By what we have seen in \S2 and \S3, the categories ${\rm Mod}_R$\,, ${\rm Topm}_R$\,, ${\rm Ltm}_R$\,, Ahtg, Ban, Lcs, ${\rm Lcs}_K$\,, ${\rm Borm}_R$\,, Cbvs and ${\rm Cbvs}_K$ are preadditive and satisfy axioms (d) and (e). In the next examples we shall see that they also satisfy axiom (f), and hence are semiabelian.

\bigskip

\noindent\textbf{Example 4.2.} ${\rm Mod}_R$ is a semiabelian category.

\medskip

Indeed, let $u \in {\rm Mor}_{{\rm Mod}_R}(A,B)$ be arbitrary. Then ${\rm Coim}(u) = {\rm Coker}({\rm Ker}(u)) = A/{\rm Ker(u)}$\,, ${\rm Im}(u) = {\rm Ker}({\rm Coker}(u)) = \{u(x); x \in A\}$ and $\overline{u}\colon {\rm Coim}(u) \to {\rm Im}(u)$ is the $R$-linear mapping given by $\overline{u}(x+{\rm Ker}(u)) = u(x)$ for $x \in A$.

\medskip

Now, let $C \in {\rm Ob}({\rm Mod}_R)$ be arbitrary, and consider the sequences
$$
0 \longrightarrow {\rm Mor}_{{\rm Mod}_R}(C,A/{\rm Ker}(u)) \overset{(\overline{u})^*}{\longrightarrow} {\rm Mor}_{{\rm Mod}_R}(C, {\rm Im}(u))
$$
and
$$
0 \longrightarrow {\rm Mor}_{{\rm Mod}_R}({\rm Im}(u),C) \overset{(\overline{u})_*}{\longrightarrow} {\rm Mor}_{{\rm Mod}_R}(A/{\rm Ker}(u), C).
$$

\smallskip

\noindent If $v \in {\rm Mor}_{{\rm Mod}_R}(C,A/{\rm Ker}(u))$ and $(\overline{u})^*(v) = \overline{u}\circ v = 0_{C\,{\rm Im}(u)}$\,, then $v = 0_{C\,A/{\rm Ker}(u)}$\,. In fact, if $t \in C$ is arbitrary, there is an $x \in A$ with $v(t) = x+{\rm Ker}(u)$; thus $\overline{u}(v(t)) = u(x) = 0$ (that is, $x \in {\rm Ker}(u))$, and $v(t) = {\rm Ker}(u)$. This shows the exactness of the first sequence. And, if $w \in {\rm Mor}_{{\rm Mod}_R}({\rm Im}(u),C)$ and $(\overline{u})_*(w) = w \circ \overline{u} = 0_{A/{\rm Ker}(u)\,C}$\,, then $w = 0_{{\rm Im}(u)\,C}$\,. In fact, if $x \in A$ is arbitrary,
$$
0 = (w \circ \overline{u})(x+{\rm Ker}(u)) = w(u(x)).
$$
This shows the exactness of the second sequence. Therefore, in view of Remark 3.3, the morphism $\overline{u}$ is bijective (more precisely, it is easily seen that $\overline{u}$ is an isomorphism in ${\rm Mod}_R)$.

\medskip

By the same procedure (and with the pertinent modifications), we conclude that the categories ${\rm Topm}_R$\,, ${\rm Ltm}_R$\,, Lcs, ${\rm Lcs}_K$\,, ${\rm Borm}_R$\,, Cbvs and ${\rm Cbvs}_K$ are semiabelian.

\bigskip

\noindent{\bf Example 4.3.} Ahtg\, is a semiabelian category.

\medskip

Indeed, let $u \in {\rm Mor}_{{\rm Ahtg}}(A,B)$ be arbitrary. Then ${\rm Coim}(u) = {\rm Coker}({\rm Ker}(u)) = A/{\rm Ker}(u)$ (endowed with the quotient topology; recall Examples 3.10 and 3.22), and ${\rm Im}(u) = {\rm Ker}({\rm Coker}(u)) = \{\overline{u(x); x \in A}\}$ (endowed with the topology induced by that $B$). Let $C \in {\rm Ob}({\rm Ahtg})$ be arbitrary, and consider the sequences
$$
0 \longrightarrow {\rm Mor}_{{\rm Ahtg}}(C,A/{\rm Ker}(u)) \overset{(\overline{u})^*}{\longrightarrow} {\rm Mor_{\rm Ahtg}}(C, {\rm Im}(u))
$$
and
$$
0 \longrightarrow {\rm Mor}_{{\rm Ahtg}}({\rm Im}(u),C) \overset{(\overline{u})_*}{\longrightarrow} {\rm Mor}_{{\rm Ahtg}}(A/{\rm Ker}(u),C),
$$
where $\overline{u}$ is the morphism in\, Ahtg\, given by $\overline{u}(x+{\rm Ker}(u)) = u(x)$ for $x \in A$. The exactness of the first sequence follows exactly as in Example 4.2. Moreover, if $w \in {\rm Mor}_{{\rm Ahtg}}({\rm Im}(u),C)$ and $(\overline{u})_*(w) = 0_{A/{\rm Ker}(u)C}$\,, then $w = 0_{{\rm Im}(u)C}$\,. In fact, if $x \in A$ is arbitrary, $w(u(x))=0$ as we have seen above, and the continuity of $w$ implies $w=0_{{\rm Im}(u)C}$\,. This shows the exactness of the second sequence. Therefore, in view of Remark 3.3, the morphism $\overline{u}$ is bijective, and (f) holds.

\bigskip

\noindent\textbf{Example 4.4.} Ban\, is a semiabelian category.

\medskip

Indeed, let $u \in {\rm Mor}_{{\rm Ban}}(A,B)$ be arbitrary. Then ${\rm Coim}(u) = A/{\rm Ker}(u)$ (endowed with the quotient norm), and ${\rm Im}(u) = \{\overline{u(x); x \in A}\}$ (endowed with the norm induced by that of $B$). By arguing exactly as in the preceding example, we conclude that the morphism $\overline{u}$ is bijective, and (f) holds.

\bigskip

\noindent\textbf{Proposition 4.5.} \textit{Let $\mathcal C$ be a semiabelian category, $u \in {\rm Mor}_{\mathcal C}(A,B)$ and $v \in {\rm Mor}_{\mathcal C}(B,C)$. Then the following conditions are equivalent:}

\smallskip

\noindent (a)\, $vu = 0_{AC}$\,;

\smallskip

\noindent (b)\, ${\rm Im}(u) \le {\rm Ker}(v)$;

\smallskip

\noindent (c)\, ${\rm Coim}(v) \le {\rm Coker}(u)$.
$$
\begin{matrix}
A \overset{j}{\longrightarrow} {\rm Coim}(u) \overset{\overline{u}}{\longrightarrow} {\rm Im}(u) \overset{i}{\longrightarrow} &B& \overset{v}{\longrightarrow} C\\
&\,\,\uparrow i'&\\
&{\rm Ker}(v)&
\end{matrix}
$$

\smallskip

\noindent\textbf{Proof.} Since (b) and (c) are dual assertions, it is enough to prove that (a) and (b) are equivalent.

\medskip

\noindent (a) $\Rightarrow$ (b):\, Since $0_{AC} = vu = v(i\,\overline{u}\,j) = (v\,i\,\overline{u})j = 0_{{\rm Coim}(u)\,C} j$ and since $j$ is surjective, we have $(vi)\overline{u} = 0_{{\rm Coim}(u)\,C} = 0_{{\rm Im}(u)\,C} \overline{u}$. Thus, by the surjectivity of $\overline{u}$, $vi = 0_{{\rm Im}(u)\,C}$\,, and Proposition 3.5 furnishes a morphism $w\colon {\rm Im}(u) \to {\rm Ker}(v)$ such that $i = i' w$; hence ${\rm Im}(u) \le {\rm Ker}(v)$.

\medskip

\noindent (b) $\Rightarrow$ (a):\, By hypothesis we can write $i = i'w$, for $w$ as above. Finally, the relations $vu = (vi)(\overline{u} j)$ and $vi = (vi')w = 0_{{\rm Ker}(v)\,C} w = 0_{{\rm Im}(u)\,C}$ give $vu = 0_{AC}$\,.

\medskip

It seems that the notion of a strict morphism was first considered by Weil
[35] in his celebrated memoir on topological groups. In our context, such a
concept reads as follows:

\bigskip

\noindent\textbf{Definition 4.6}\,[26]. Let $\mathcal C$ be a semiabelian category. A morphism $u$ in $\mathcal C$ is said to be \textit{strict\/} if the corresponding morphism $\overline{u}$ is an isomorphism.

\medskip

It is easily seen that, if $u \in {\rm Mor}_{\mathcal C}(A,B)$ is arbitrary, then ${\rm Ker}(u) \overset{i}{\longrightarrow} A$ and $B \overset{j}{\longrightarrow} {\rm Coker}(u)$ are strict morphisms.

\medskip

For the rest of our work we shall use the following convention:\, if $\mathcal C$ is a semiabelian category, $A \in {\rm Ob}(\mathcal C)$ and $(A',i')$, $(A'',i'')$ are subobjects of $A$ (resp. $(Q',p')$, $(Q'',p'')$ are quotients of $A$), the symbol $A'=A''$ (resp. $Q'=Q''$) will mean that $A'$ and $A''$ (resp. $Q'$ and $Q''$) are isomorphic.

\bigskip

\noindent\textbf{Proposition 4.7.} \textit{Let $\mathcal C$ be a semiabelian category and $A \in {\rm Ob}(\mathcal C)$. Then the following assertions hold:}
\begin{itemize}
\item[{(a)}]\, \textit{if $(A',i')$ is a strict subobject of $A$ and $({\rm Coker}(i'),j')$ is the cokernel of $i'$, then ${\rm Ker}(j') = A'$;}
\item[{(b)}]\, \textit{if $(Q',p')$ is a strict quotient of $A$ and $({\rm Ker}(p'),i')$ is the kernel of $p'$, then ${\rm Coker}(i') = Q'$.}
\end{itemize}

\medskip

\noindent\textbf{Proof.} (a):\, Let us consider the sequence
$$
A' \overset{\ell}{\longrightarrow} {\rm Coim}(i') \overset{\overline{i'}}{\longrightarrow} {\rm Im}(i').
$$
By Proposition 3.30(b), $\ell$ is an isomorphism; and, by hypothesis, $\overline{i'}$ is an isomorphism. Thus the morphism $\overline{i'} \ell$ is an isomorphism. Consequently,
$$
A' = {\rm Im}(i') = {\rm Ker}({\rm Coker}(i')) = {\rm Ker}(j').
$$

\smallskip

\noindent (b):\, Follows from (a), by duality.

\bigskip

\noindent\textbf{Definition 4.8.} Let $\mathcal C$ be a semiabelian category and $A \in {\rm Ob}(\mathcal C)$. For each strict subobject $(A',i')$ of $A$, we shall write ${\rm Coker}(i') = A/A'$.

\medskip

The proof of the next result follows the lines of that of Theorem 2.13 of [12].

\bigskip

\noindent\textbf{Proposition 4.9.} \textit{Let $\mathcal C$ be a semiabelian category and $A \in {\rm Ob}(\mathcal C)$. Then any pair $(A',i')$, $(A'',i'')$ of strict subobjects of $A$ admits an infimum $($denoted by $A' \cap A'')$ in the preordered class of strict subobjects of $A$.}

\medskip

\noindent\textbf{Proof.} Put ${\rm Coker}(i') = (A/A', j')$ and
$v = j'i''\colon A'' \to A/A'$\,, and let $({\rm Ker}(v),i)$ be the kernel
of $v$. Then $({\rm Ker}(v),i'' i)$ is a strict subobject of $A$ by Lemma 6
of [18], and ${\rm Ker}(v) \le A''$. We claim that ${\rm Ker}(v) \le A'$. In fact, since
$$
0_{{\rm Ker}(v)A/A'} = vi = (j' i'')i = j'(i'' i)
$$
and since ${\rm Ker}(j')=A'$ by Proposition 4.7(a), there is a morphism $w\colon{\rm Ker}(v) \to A'$ in $\mathcal C$ such that the diagram
$$
\begin{matrix}
{\rm Ker}(v)  &\overset{w}{\longrightarrow}& A'\\
i \downarrow & & \downarrow i'\\
A'' &\underset{i''}{\longrightarrow}& A
\end{matrix}  
$$
is commutative by Proposition 3.5. Thus ${\rm Ker}(v) \le A'$.

\medskip

Now, let $(X,k)$ be a strict subobject of $A$ such that $X \le A'$ and $X \le A''$. We claim that $X \le {\rm Ker}(v)$. Indeed, since $X \le A'$, there is a morphism $\theta_1\colon X \to A'$ in $\mathcal C$ making the diagram
$$
\begin{matrix}
X  &\overset{\theta_1}{\longrightarrow}& A'\\
k \searrow & & \swarrow i'\\
&A&
\end{matrix}  
$$
commutative. And, since $X \le A''$, there is a morphism $\theta_2\colon X \to A''$ in $\mathcal C$ making the diagram
$$
\begin{matrix}
X  &\overset{\theta_2}{\longrightarrow}& A''\\
k \searrow & & \swarrow i''\\
&A&
\end{matrix}  
$$
commutative. On the other hand, since
$$
v\theta_2 = (j' i'')\theta_2 = j'(i'' \theta_2) = j' k = j'(i' \theta_1) = (j' i')\theta_1 = 0_{A'A/A'}\,\theta_1 = 0_{XA/A'}
$$
and since $({\rm Ker}(v),i)$ is the kernel of $v$, Proposition 3.5 guarantees the existence of morphism $t\colon X \to {\rm Ker}(v)$ in $\mathcal C$ making the diagram
$$
\begin{matrix}
X  &\overset{t}{\longrightarrow}& {\rm Ker}(v)\\
\theta_2 \searrow & & \swarrow i\\
&A''&
\end{matrix}  
$$
commutative. Consequently,
$$
k = i'' \theta_2 = i''(it) = (i'' i)t,
$$
proving that $X \le {\rm Ker}(v)$. This completes the proof. 

\bigskip

\noindent\textbf{Corollary 4.10.} \textit{Let $\mathcal C$ be a semiabelian category and $A \in {\rm Ob}(\mathcal C)$. Then the following assertions hold:}
\begin{itemize}
\item[{(a)}] \textit{every pair $(Q',p')$, $(Q'',p'')$ of strict quotients of $A$ admits an infimum $($denoted by $Q' \cap Q'')$ and a supremum $($denoted by $Q' \cup Q'')$ in the preordered class of strict quotients of $A$;}
\item[{(b)}] \textit{every pair $(A',i')$, $(A'',i'')$ of strict subobjects of $A$ admits a supremum $($denoted by $A' \cup A'')$ in the preordered class of strict subobjects of $A$.}
\end{itemize}

Now we can state the following 

\bigskip

\noindent\textbf{Theorem 4.11.} \textit{Let $\mathcal C$ be a semiabelian category and $A \in {\rm Ob}(\mathcal C)$, and let $(A',i')$, $(A'',i'')$ be strict subobjects of $A$. Then there exists a bijective morphism}
$$
A''/(A' \cap A'') \longrightarrow (A' \cup A'')/A'
$$
\textit{in $\mathcal C$.}

\medskip

In order to prove Theorem 4.11 we shall need two auxiliary results.

\bigskip

\noindent\textbf{Lemma 4.12.} \textit{Let $\mathcal C$ be a semiabelian category and let $u \in {\rm Mor}_{\mathcal C}(A,B)$ be such that ${\rm Ker}(u)=A$. Then $u = 0_{AB}$\,.}

\medskip

\noindent\textbf{Proof.} We have ${\rm Ker}(u) \overset{i}{\longrightarrow} A \overset{u}{\longrightarrow} B$, with $ui = 0_{{\rm Ker}(u)B}$\,, $i$ being an isomorphism by hypothesis. Therefore
$$
0_{AB} = (ui)i^{-1} = u(ii^{-1}) = u.
$$

\medskip

\noindent\textbf{Lemma 4.13.} \textit{Let $\mathcal C$, $A$, $(A',i')$ and $(A'', i'')$ be as in Theorem $4.11$, and consider the sequence}
$$
A'' \overset{k}{\longrightarrow} A' \cup A'' \overset{\ell}{\longrightarrow} (A' \cup A'')/A'.
$$
\textit{Then the morphism $u = \ell k$ is surjective.}

\medskip

\noindent\textbf{Proof.} Let $w\colon (A' \cup A'')/A' \to C$ be a morphism in $\mathcal C$ such that $u_*(w) = wu = (w\ell)k = 0_{A''C}$\,. We have to show that $w = 0_{(A'\cup A'')/A'\,C}$ (Remark 3.3). But, since $\ell$ is surjective, it is enough to show that $w\ell = 0_{A'\cup A''\,C}$\,. Indeed, by Proposition 4.5, ${\rm Im}(k) \le {\rm Ker}(w\ell)$. On the other hand, by Proposition 4.7(a),
$$
{\rm Im}(k) = {\rm Ker}({\rm Coker}(k)) = A'';
$$
thus $A'' \le {\rm Ker}(w\ell)$. Moreover, by Propositions 4.7(a) and 3.32(a), $A' = {\rm Ker}(\ell) \le {\rm Ker}(w\ell)$. Consequently, $A' \cup A'' \le {\rm Ker}(w\ell)$. But, since ${\rm Ker}(w\ell) \le A' \cup A''$, we get ${\rm Ker}(w\ell) = A' \cup A''$. Therefore, by Lemma 4.12, $w\ell = 0_{A'\cup A''\,C}$\,.

\medskip

Now, let us turn to the

\bigskip

\noindent\textbf{Proof of Theorem 4.11.} We may assume that $A = A' \cup A''$. Let $u$ be as in the proof of Lemma 4.13. As we have seen in the proof of Proposition 4.9, ${\rm Ker}(u) = A' \cap A''$; thus
$$
{\rm Coim}(u) = {\rm Coker}({\rm Ker}(u)) = A''/(A'\cap A'')\,.
$$
On the other hand, the morphism
$$
A''/(A'\cap A'') \overset{\overline{u}}{\longrightarrow} {\rm Im}(u)
$$
is bijective. Finally, by Lemma 4.13 and Proposition 3.30(d), ${\rm Im}(u) = (A' \cup A'')/A'$, which concludes the proof.

\bigskip

\noindent\textbf{Remark 4.14.} The bijective morphism
$A''/(A'\cap A'') \to (A' \cup A'')/A'$
considered in the proof of Theorem 4.11 is not necessarily an isomorphism;
for instance, consider the semiabelian category of abelian topological groups
and the example given in [4], p.27.

\medskip

In order to prove a version of the first isomorphism theorem for semiabelian
categories (see also [26]) we shall need the following

\bigskip

\noindent\textbf{Proposition 4.15.} \textit{Let $\mathcal C$ be a semiabelian category and $A \in {\rm Ob}(\mathcal C)$, and let $(A',i')$, $(A'',i'')$ be strict subobjects of $A$ such that $A' \le A''$. Then there exists a unique strict surjective morphism $v\colon A/A' \to A/A''$ in $\mathcal C$ making the diagram}
$$
\begin{matrix}
&A&\\
&j' \swarrow \qquad \searrow j''&\\
&A/A' \underset{v}{\longrightarrow} A/A''&
\end{matrix}
$$
\textit{commutative, where $j' = {\rm Coker}(i')$ and $j'' = {\rm Coker}(i'')$.}

\medskip

\noindent\textbf{Proof.} By hypothesis, there is a morphism $w\colon A' \to A''$ in $\mathcal C$ such that $i' = i'' w$ and, by Proposition 4.7(a), ${\rm Ker}(j') = A'$ and ${\rm Ker}(j'') = A''$. On the other hand, since
$$
j''(i'' w) = (j'' i'')w = 0_{A''\, A/A''}\,w = 0_{A'\,A/A''}\,,
$$
Proposition 3.17 furnishes a morphism $v\colon A/A' \to A/A''$ in
$\mathcal C$ such that $j'' = vj'$. Since the surjectivity of $v$ follows
from assertion (b) after Definition 2.15 and the strictness of $v$ follows
from Lemma 4 of [18], the existence is proved.

\medskip

To prove the uniqueness, let $v,v' \in {\rm Mor}_{\mathcal C}(A/A', A/A'')$ be such that $j'' = vj'$ and $j'' = v'j'$. Then $vj' = v'j'$, and the surjectivity of $j'$ gives $v=v'$.

\bigskip

\noindent\textbf{Theorem 4.16.} \textit{Let $\mathcal C$ be a semiabelian category and $A \in {\rm Ob}(\mathcal C)$. Let $(A'',i'')$ be a strict subobject of $A$ and $(A',i)$ a strict subobject of $A''$, and let $v\colon A/A' \to A/A''$ be as in Proposition $4.15$. Then there exists a sequence}
$$
A''/A' \longrightarrow {\rm Ker}(v) \overset{\ell}{\longrightarrow} A/A' \overset{v}{\longrightarrow} A/A''
$$
\textit{of morphisms in $\mathcal C$, the morphism on the left being bijective. Moreover, $(A/A')\big/{\rm Ker}(v)$ and $A/A''$ are isomorphic.}

\medskip

\noindent\textbf{Proof.} Put $i' = i''i$. By Lemma 6 of [18], $(A',i')$ is a strict subobject of $A$, and hence $A/A'$ makes sense. Let $\alpha$ be the morphism $A'' \overset{i''}{\longrightarrow} A \overset{j'}{\longrightarrow} A/A'$ and consider the sequence
$$
A' \overset{i}{\longrightarrow} A'' \overset{\alpha}{\longrightarrow} A/A' \overset{v}{\longrightarrow} A/A''
$$
of morphisms in $\mathcal C$. We claim that ${\rm Im}(i) = {\rm Ker}(\alpha)$ and ${\rm Im}(\alpha) = {\rm Ker}(v)$. Indeed, by Proposition 4.7(a), ${\rm Im}(i) = {\rm Ker}({\rm Coker}(i)) = A'$. And, since
$$
\alpha i = (j'i'')i = j'(i'' i) = j'i' = 0_{A'\,A/A'}\,,
$$
Proposition 4.5 furnishes ${\rm Im}(i) \le {\rm Ker}(\alpha)$. Let $s$ be the injective morphism ${\rm Ker}(\alpha) \longrightarrow A'' \overset{i''}{\longrightarrow} A$. Since the morphism
$$
{\rm Ker}(\alpha ) \longrightarrow \underbrace{A'' \longrightarrow A \longrightarrow A/A'}_{\alpha}
$$
is equal to $0_{{\rm Ker}(\alpha)A/A'}$\,, that is, since the morphism
$$
{\rm Ker}(\alpha) \overset{s}{\longrightarrow} A \overset{j'}{\longrightarrow} A/A'
$$
is equal to $0_{{\rm Ker}(\alpha)A/A'}$\,, Proposition 3.5 implies the existence of a morphism $w\colon {\rm Ker}(\alpha) \to A'$ in $\mathcal C$ making the diagram
$$
\begin{matrix}
&{\rm{Ker}(\alpha)} \overset{w}{\longrightarrow} A'&\\
&s \searrow \qquad \swarrow i'&\\
&A&
\end{matrix}
$$
commutative. Thus ${\rm Ker}(\alpha) \le A' = {\rm Im}(i)$, and hence ${\rm Im}(i) = {\rm Ker}(\alpha)$. Therefore, by duality, ${\rm Im}(\alpha) = {\rm Ker}(v)$. Moreover, by hypothesis, the morphism $\overline{\alpha}\colon {\rm Coim}(\alpha) \to {\rm Im}(\alpha)$ is bijective, where
$$
{\rm Coim}(\alpha) = {\rm Coker}({\rm Ker}(\alpha)) = A''/A'.
$$
Finally, ${\rm Coker}(\ell) = (A/A')\big/{\rm Ker}(v)$ is isomorphic to $A/A''$ in view of Proposition 4.7(b).

This completes the proof.

\medskip

Before proceeding we would like to mention that the proofs of Propositions 4.5 and 4.7 and of Theorem 4.16 are essentially contained in [2].

\bigskip

\noindent\textbf{Remark 4.17.} (a)\, Let $A \in {\rm Ob}({{\rm Mod}_R})$ and let $A' \subset A''$ be a submodules of $A$ (resp.\linebreak 
\!\!\!\!\!
\noindent $A \in {\rm Ob}{({\rm Topm}_R)}$ and let $A' \subset A''$ be submodules of $A$ endowed with the induced topology,

\noindent $A \in {\rm Ob}{({\rm Ltm}_R)}$ and let $A' \subset A''$ be submodules of $A$ endowed with the induced topology,

\noindent $A \in {\rm Ob}{({\rm Lcs})}$ and let $A' \subset A''$ be subspaces of $A$ endowed with the induced topology,

\noindent $A \in {\rm Ob}{({\rm Lcs}_K)}$ and let $A' \subset A''$ be subspaces of $A$ endowed with the induced topology,

\noindent $A \in {\rm Ob}{({\rm Borm}_R)}$ and let $A' \subset A''$ be submodules of $A$ endowed with the induced bornology,

\noindent $A \in {\rm Ob}{({\rm Cbvs})}$ and let $A' \subset A''$ be subspaces of $A$ endowed with the induced bornology,

\noindent $A \in {\rm Ob}{({\rm Cbvs}_K})$ and let $A' \subset A''$ be subspaces of $A$ endowed with the induced bornology).

\noindent If $v\colon A/A' \to A/A''$ is as in Theorem 4.16, it is easily seen that ${\rm Ker}(v) = A''/A'$, and therefore $(A/A')/(A''/A')$ is isomorphic to $A/A''$.

\smallskip

\noindent (b)\, Let $A \in {\rm Ob}({\rm Ahtg})$ and let $A' \subset A''$ be closed subgroups of $A$ endowed with the induced topology (resp. $A \in {\rm Ob}({\rm Ban})$ and let $A' \subset A''$ be closed subspaces of $A$ endowed with the induced norm). 
If $v\colon A/A' \to A/A''$ is as in Theorem 4.16, it is easily seen that ${\rm Ker}(v) = A''/A'$, and therefore $(A/A')/(A''/A')$ is isomorphic to $A/A''$.

\medskip

We shall conclude our paper with a few comments on the notion of an abelian
category [12], introduced by Buchsbaum [7] and Grothendieck [14].

\bigskip

\noindent\textbf{Definition 4.18.} A category $\mathcal C$ is said to be \textit{abelian\/} if it is preadditive, satisfies axioms (d) and (e) and the following axiom:

\smallskip

\noindent (f1)\, for all $A,B \in {\rm Ob}(\mathcal C)$ and for all $u \in {\rm Mor}_{\mathcal C}(A,B)$, the morphism $\overline{u}$ is an isomorphism.

\medskip

Since every isomorphism is bijective, it follows that every abelian category is semiabelian.

\bigskip

\noindent\textbf{Remark 4.19.} In view of Propositions 3.30(b),(d) and 3.31 and the fact that every isomorphism is bijective, axiom (f1) is equivalent to the following axiom:

\smallskip

\noindent (f2)\,(a)\, for every morphism $u$ in $\mathcal C$, the morphism $\overline{u}$ is bijective; (b)\, every bijective morphism in $\mathcal C$ is an isomorphism.

\medskip 

Let us mention a few examples of abelian categories.

\bigskip

\noindent\textbf{Example 4.20.} The category ${\rm Mod}_R$ is abelian.

\bigskip

\noindent\textbf{Example 4.21.} The category of finite abelian groups is abelian.

\bigskip

\noindent\textbf{Example 4.22.} The category of finite abelian $p$-groups [20] is abelian.

\bigskip

\noindent\textbf{Example 4.23.} The category of vector bundles [19] is abelian.

\bigskip

\noindent\textbf{Example 4.24.} The category of sheaves of abelian groups
over a topological space [15] is abelian.

\medskip

Now, let us see some examples of semiabelian categories which are not abelian.

\bigskip

\noindent\textbf{Example 4.25.} Let $K$ be a discrete field. Then the semiabelian category of topological vector spaces over $K$ is not abelian.

\medskip

In fact, let $A$ be a vector space over $K$, with $A \ne \{0\}$, and let $\tau_1$ (resp. $\tau_2$) be the discrete topology (resp. the chaotic topology) on $A$. Then $x \in (A,\tau_1) \mapsto x \in (A,\tau_2)$ is a bijective morphism in the category of topological vector spaces over $K$ which is not an isomorphism. Thus, by Remark 4.19, this category is not abelian.

\bigskip

\noindent\textbf{Example 4.26.} Let $K$ be a discrete field. Then the semiabelian category of linearly topologized spaces over $K$ is not abelian.

\medskip

In fact, it suffices to argue as in the preceding example, by observing that $(A,\tau_1)$ and $(A,\tau_2)$ are linearly topologized spaces over $K$.

\bigskip

\noindent\textbf{Example 4.27.} The semiabelian category Ahtg is not abelian.

\medskip

In fact, it suffices to recall Example 2.16 and Remark 4.19.

\bigskip

\noindent\textbf{Example 4.28.} The semiabelian category Ban is not abelian.

\medskip

In fact, consider the Banach spaces $A = (\ell^1, || \cdot ||_1)$ and $B = (c_0, || \cdot ||_0)$, and let $u \colon A \to B$ be the linear mapping given by $u(x)=x$ for $x \in A$. Since $||u(x)||_0 \le ||x||_1$ for all $x \in A$, \, $u \in {\rm Mor}_{\rm Ban}(A,B)$. Moreover, since $A$ is dense in $B$, it follows that $u$ is a bijective morphism in Ban. Finally, $u$ is not an isomorphism, and Remark 4.19 implies that Ban is not abelian.

\bigskip

\noindent\textbf{Example 4.29.} The semiabelian category Lcs is not abelian.

\medskip

In fact, let $A$ be an infinite-dimensional normed space endowed with the locally convex topology $\tau$ defined by its norm, and let $\sigma(A,A')$ be the weak topology on $A$. Then $x \in (A,\tau) \mapsto x \in (A, \sigma(A,A'))$ is a bijective morphism in Lcs which is not an isomorphism. Thus, by Remark 4.19,\, Lcs is not abelian.

\bigskip

\noindent\textbf{Example 4.30.} The semiabelian category ${\rm Lcs}_K$ is not abelian.

\medskip

In fact, consider the vector space $A = K^{(\mathbb N)}$ over $K$ endowed
with the following locally $K$-convex topologies $\tau_1$ and $\tau_2$:
$\tau_1$ is the direct sum topology [33, p.\ 268] (which is not metrizable
by Theorem 3.13 of [33]) and $\tau_2$ is the topology on $A$ induced by the product topology on $K^{\mathbb N}$ (which is obviously metrizable). Then $x \in (A,\tau_1) \mapsto x \in (A,\tau_2)$ is a bijective morphism in ${\rm Lcs}_K$ which is not an isomorphism. Therefore, by Remark 4.19, \, ${\rm Lcs}_K$ is not abelian.

\bigskip

\noindent\textbf{Example 4.31.} Let $K$ be a discrete field. Then the semiabelian category of bornological vector spaces over $K$ is not abelian.

\medskip

In fact, let $A$ be an infinite-dimensional vector space over $K$, and let $\mathcal{B}_1$ (resp. $\mathcal{B}_2$) be the vector bornology having the finite-dimensional subspaces of $A$ as a fundamental system of bounded sets (resp. the vector bornology consisting of all subsets of $A$). Then $x \in (A,\mathcal{B}_1) \mapsto x \in (A,\mathcal{B}_2)$ is a bijective morphism in the category of bornological vector spaces over $K$ which is not an isomorphism. Therefore, by Remark 4.19, this category is not abelian.

\bigskip

\noindent\textbf{Example 4.32.} The semiabelian category ${\rm Cbvs}_K$ is not abelian if $K$ is spherically complete.

\medskip

In fact, consider $(A,\tau_1)$ and $(A,\tau_2)$ as in Example 4.30 and let
$\mathcal{B}_1$ (resp. $\mathcal{B}_2$) be the $K$-convex bornology
consisting of all $\tau_1$-bounded (resp. $\tau_2$-bounded) subsets of $A$.
Then $x \in (A,\mathcal{B}_1) \mapsto x \in (A,\mathcal{B}_2)$ is a bijective
morphism in ${\rm Cbvs}_K$ which is not an isomorphism (for the boundedness
of $x \in (A,\mathcal{B}_2) \mapsto x \in (A,\mathcal{B}_1)$ would imply the
continuity of\linebreak $x \in (A,\tau_2) \mapsto x \in (A,\tau_1)$ in view
of the remark after Theorem 4.30 of [33]). Therefore, by Remark 4.19,\, ${\rm Cbvs}_K$ is not abelian.

\medskip

In view of [9] we can argue as in Examples 4.30 and 4.32 to conclude:

\bigskip

\noindent\textbf{Example 4.33.} The semiabelian category Cbvs is not abelian.

\medskip

We close our work by deriving the well-known isomorphism theorems for abelian
categories [2,8,12].

\bigskip

\noindent\textbf{Theorem 4.34.} \textit{Let $\mathcal C$ be an abelian category and $A \in {\rm Ob}(\mathcal C)$, and let $(A',i')$, $(A'',i'')$ be subobjects of $A$. Then}
$$
A''\big/(A' \cap A'') \quad\rm{and}\quad (A' \cup A'')\big/A'
$$
\textit{are isomorphic.}

\medskip

\noindent\textbf{Proof.} Follows immediately from Theorem 4.11 and Remark 4.19.

\bigskip

\noindent\textbf{Theorem 4.35.} \textit{Let $\mathcal C$ be an abelian category and $A \in {\rm Ob}(\mathcal C)$. If $(A'',i'')$ is a subobject of $A$ and $(A',i)$ is a subobject of $A''$, then}
$$
(A/A')\big/(A''/A') \quad\rm{and}\quad A/A''
$$
\textit{are isomorphic.}

\medskip

\noindent\textbf{Proof.} Follows immediately from Theorem 4.16 and Remark 4.19.

\newpage

\centerline{\large\bf References}

\vglue .1in

\begin{itemize}
\item[{[1]}] M. Akkar,\ Espaces vectoriels bornologiques $K$-convexes, \ Indag. Math. 73 (1970), 82-95.
\item[{[2]}] A. Andreotti,\ Généralités sur les catégories abéliennes, \ Séminaire A. Grothendieck (Algèbre Homologique), Exposé 2, Faculté des Sciences de Paris (1957).
\item[{[3]}] J. Bonet \& S. Dierof,
             The pullback for bornological and ultrabornological spaces,
             Note Mat.\ 25 (2006), 63-67.
\item[{[4]}] N. Bourbaki,\ Topologie Générale, Chapitres 3 et 4, Troisième édition, \ Actualités Scientifiques et Industrielles 1084, Hermann (1960).
\item[{[5]}] N. Bourbaki,\ Algèbre, Chapitre 2, Troisième édition, \ Actualités Scientifiques et Industrielles 1236, Hermann (1961).
\item[{[6]}] N. Bourbaki,\ Topologie Générale, Chapitres 1 et 2, Quatrième édition, \ Actualités Scientifiques et Industrielles 1142, Hermann (1965).
\item[{[7]}] D.A. Buchsbaum,\ Exact categories and duality, \ Trans. Amer. Math. Soc. 80 (1955), 1-34.
\item[{[8]}] I. Bucur \& A. Deleanu,\ Introduction to the Theory of Categories and Functors, \ Texts, Monographs and Tracts in Pure and Applied Mathematics XIX, John Wiley and Sons (1974).
\item[{[9]}] J. Dieudonné \& L. Schwartz,\ La dualité dans les espaces $(\mathcal L)$ et $(\mathcal{L}\mathcal{F})$, \ Ann. Inst. Fourier 1 (1949), 61-101.
\item[{[10]}] J. Dieudonné,\ Foundations of Modern Analysis, \ Second (enlarged and corrected) printing, Academic Press (1969).
\item[{[11]}] S. Eilenberg \& S. Mac Lane,\ General theory of natural equivalences, \ Trans. Amer. Math. Soc. 58 (1945), 231-294.
\item[{[12]}] P. Freyd,\ Abelian Categories, \ Harper and Row (1964); republished in: Reprints in Theory and Applications of Categories 3 (2003), 23-164.
\item[{[13]}] P. Gabriel,\ Des catégories abéliennes, \ Bull. Soc. Math. France 90 (1962), 323-448.
\item[{[14]}] A. Grothendieck,\ Sur quelques points d'algèbre homologique, \ Tôhoku Math. J.\linebreak 9 (1957), 119-221.
\item[{[15]}] R. Hartshorne,\ Algebraic Geometry, \ Graduate Texts in Mathematics 52, Springer-Verlag (1977).
\item[{[16]}] H. Hogbe-Nlend,\ Bornologies and Functional Analysis, \ Notas de Matemática 62, North-Holland (1977).
\item[{[17]}] G. K\"{o}the,\ Topological Vector Spaces I, \ Grundlehren der mathematischen Wissenschaften 159, Springer-Verlag (1969).
\item[{[18]}] V.I. Kuz'minov \& A. Yu. Cherevikin,\ Semiabelian categories, \ Siberian Math. J. 13 (1972), 895-902.
\item[{[19]}] S. Lang,\ Introduction to Differentiable Manifolds, \ Interscience (1962).
\item[{[20]}] S. Lang,\ Algebra, \ Revised third edition, Graduate Texts in Mathematics 211, Springer-Verlag (2002).
\item[{[21]}] S. Lefschetz,\ Algebraic Topology, \ Amer. Math. Soc. Colloq. Publ. XXVII (1942).
\item[{[22]}] V.P. Palamodov,\ Homological methods in the theory of locally convex spaces, \ Russian Math. Surveys 26 (1971), 1-64.
\item[{[23]}] D.P. Pombo Jr.,\ Linear topologies and linear bornologies on modules, \ J. Indian Math. Soc. 59 (1993), 107-117; errata ibid. 60 (1994), 191.
\item[{[24]}] D.P. Pombo Jr.,\ Module bornologies, \ Comment. Math. (Prace Mat.) 36 (1996), 179-188.
\item[{[25]}] D.P. Pombo Jr.,\ On the notion of a semiabelian category (2013), \ preprint.
\item[{[26]}] A. Ra$\overset{\vee}{\,\!\rm \i}\!$kov,\ Semiabelian categories, \ Soviet Math. J. 10 (1969), 1242-1245.
\item[{[27]}] A.P. Robertson \& W. Robertson,\ Topological Vector Spaces, \ Second edition, Cambridge University Press (1973).
\item[{[28]}] W. Rump,
              A counterexample to Ra$\overset{\vee}{\,\!\rm \i}\!$kov's
              conjecture,
              Bull.\ London Math.\ Soc.\ 40 (2008), 985-994.
\item[{[29]}] W. Rump,
              Analysis of a problem of Ra$\overset{\vee}{\,\!\rm \i}\!$kov
              with applications to barreled and bornological spaces,
              J. Pure Appl.\ Algebra 215 (2011), 44-52.
\item[{[30]}] P. Schneider,\ Nonarchimedean Functional Analysis, \ Springer Monographs in Mathematics, Springer-Verlag (2002).
\item[{[31]}] J.-P. Schneiders,\ Quasi-abelian categories and sheaves, \ Mém. Soc. Math. France 76 (1999).
\item[{[32]}] Séminaire Banach, \ Lecture Notes in Mathematics 277, Springer-Verlag (1972).
\item[{[33]}] J. van Tiel,\ Espaces localement $K$-convexes, \ Indag. Math. 68 (1965), 249-289.
\item[{[34]}] S. Warner,\ Topological Fields, \ Notas de Matemática 126, North-Holland (1989).
\item[{[35]}] A. Weil,\ L'int\'egration dans les groupes topologiques et ses applications, \ Actualités Scientifiques et Industrielles 869, Hermann (1940).
\item[{[36]}] J. Wengenroth,
              The Ra$\overset{\vee}{\,\!\rm \i}\!$kov conjecture fails
              for simple analytical reasons,
              J. Pure Appl.\ Algebra 216 (2012), 1700-1703.
\end{itemize}

\end{document}